\documentclass[12pt]{article}

\usepackage[english]{babel}
\usepackage{tgtermes}
\usepackage{MnSymbol}

\usepackage{amsmath}
\usepackage{amsthm}
\usepackage{amsfonts}

\usepackage{pgfplots}
\pgfplotsset{compat=1.18}
\usepackage[hyphens]{url}
\usepackage{hyperref}

\usepackage{caption}
\usepackage{subcaption}

\usepackage{epigraph}

\usepackage{tikz}
\usetikzlibrary{arrows.meta}

\newtheorem{theorem}{Theorem}[section]

\newcommand{\Li}{\textrm{Li}_2}

\begin{document}
	\begin{center}
		\Large With what probability does an inscribed triangle contain a given point?
	\end{center}

	\begin{center}
		Abdulamin Ismailov\footnote{Higher School of Economics, \href{mailto:}{\nolinkurl{arismailov@edu.hse.ru}}}
	\end{center}

	\begin{abstract}
		Three points uniformly selected on the unit circle form a triangle containing a point $X$ at distance $r \in [0; 1]$ from its center with probability $P(r) = \frac{1}{4} - \frac{3}{2 \pi^2}\Li(r^2)$, where $\Li$ is the dilogarithm function (Jeremy Tan Jie Rui, 2018). In this paper we present an alternative proof of this fact. We also discuss a couple of other geometric probability problems where the dilogarithm function arises.
	\end{abstract}


	\section{Introduction}
	\vspace{-1em}
	\epigraph{triangle = integral}{}
	\vspace{-1em}

	With what probability do three points uniformly selected on a circle form a triangle containing its center? We see below that the triangles in question are the non-obtuse triangles.

	\begin{figure}[h]
		\centering
		\label{center}
		\begin{subfigure}{0.49\textwidth}
			\centering

			\begin{tikzpicture}
				\draw[fill] (0, 0) circle (1pt) node[above] (O) {$O$};

				\draw (0, 0) circle (60pt);
				\draw[fill] (20: 60pt) circle (1pt) node[right] (A) {$A$} node[below=6pt, left] {$\alpha$};
				\draw[fill] (160: 60pt) circle (1pt) node[left] (B) {$B$} node[below=7pt, right=15pt] {$\beta$};
				\draw[fill] (-50: 60pt) circle (1pt) node[right] (C) {$C$};
				\draw (20: 60pt) -- (160: 60pt) -- (-50: 60pt) -- (20: 60pt);
			\end{tikzpicture}
			\caption{Acute triangle}
			\label{center:acute}
		\end{subfigure}
		\hfill
		\begin{subfigure}{0.49\textwidth}
			\centering

			\begin{tikzpicture}
				\draw[fill] (0, 0) circle (1pt) node[below] (O) {$O$};

				\draw (0, 0) circle (60pt);
				\draw[fill] (40: 60pt) circle (1pt) node[right] (A) {$A$} node[below=4pt, left=10pt] {$\alpha$};
				\draw[fill] (140: 60pt) circle (1pt) node[above] (B) {$B$} node[below=8pt, right=-4pt] {$\beta$};
				\draw[fill] (-170: 60pt) circle (1pt) node[left] (C) {$C$} node[above=11pt, right=2pt] {$\gamma$};
				\draw (40: 60pt) -- (140: 60pt) -- (-170: 60pt) -- (40: 60pt);
			\end{tikzpicture}
			\caption{Obtuse triangle}
			\label{center:obtuse}
		\end{subfigure}
		\caption{}
	\end{figure}
	The problem can then be reduced to a system
	\begin{equation}
		\begin{cases}
			\alpha + \beta + \gamma = \pi \\
			\alpha \leq \frac{\pi}{2} \\ 
			\beta \leq \frac{\pi}{2} \\
			\gamma \leq \frac{\pi}{2}
		\end{cases},
	\end{equation}
	which can be visualized as follows:
	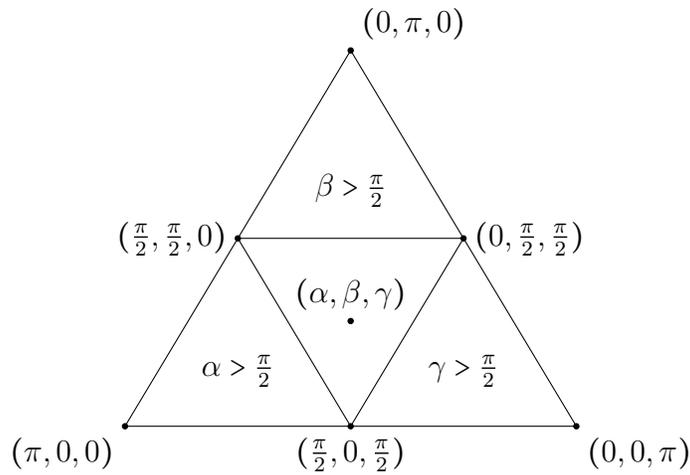
\begin{figure}[h]
		\centering
		\label{barycentric}

		\begin{tikzpicture}
			\draw[fill] (3, 0) circle (1pt) node[below right] {$(0, 0, \pi)$};
			\draw[fill] (-3, 0) circle (1pt) node[below left] {$(\pi, 0, 0)$};
			\draw[fill] (0, 5) circle (1pt) node[above right] {$(0, \pi, 0)$};

			\draw[fill] (0, 0) circle (1pt) node[below] {$(\frac{\pi}{2}, 0, \frac{\pi}{2})$};
			\draw[fill] (1.5, 2.5) circle (1pt) node[right] {$(0, \frac{\pi}{2}, \frac{\pi}{2})$};
			\draw[fill] (-1.5, 2.5) circle (1pt) node[left] {$(\frac{\pi}{2}, \frac{\pi}{2}, 0)$};

			\draw (-3, 0) -- (3, 0) -- (0, 5) -- (-3, 0);
			\draw (1.5, 2.5) -- (-1.5, 2.5) -- (0, 0) -- (1.5, 2.5);
			\draw[fill] (0, 1.4) circle (1pt) node[above] {$(\alpha, \beta, \gamma)$};

			\draw[fill] (-1.5, 1.1) circle (0pt) node[below] {$\alpha > \frac{\pi}{2}$};
			\draw[fill] (1.5, 1.1) circle (0pt) node[below] {$\gamma > \frac{\pi}{2}$};
			\draw[fill] (0, 3.5) circle (0pt) node[below] {$\beta > \frac{\pi}{2}$};
		\end{tikzpicture}

		\caption{Think of $(\alpha, \beta, \gamma)$ as of barycentric coordinates}
	\end{figure}
	We conclude that the answer is $\frac{1}{4}$. But what if instead of the center $O$ we pick an arbitrary point $X$ inside of our circle? Clearly, the answer only depends on the ratio $r = \frac{OX}{R}$, where $R$ is the circle's radius. Denote the probability in question by $P(r)$. So far we know
	\begin{equation}\label{end_values} P(0) = \frac{1}{4}, \quad P(1) = 0, \end{equation}
	the second equality corresponds to a point lying on the circle.

	Jeremy Tan Jie Rui has shown\footnote{in a Mathematics Stack Exchange post \url{https://math.stackexchange.com/questions/2641793/probability-that-a-random-triangle-with-vertices-on-a-circle-contains-an-arbitra}}
	\begin{theorem}
		\label{main_theorem}
		Three points uniformly selected on a circle with center $O$ and radius $R$ form a triangle containing a point $X$ inside it with probability
		\begin{equation}\label{main_eq} P(r) = \frac{1}{4} - \frac{3}{2 \pi^2} \Li(r^2), \end{equation}
		where $r = \frac{OX}{R}$.
	\end{theorem}

	We remind the reader that the dilogarithm function for $|z| \leq 1$ is given by power series
	\begin{equation}\label{li_2_series} \Li(z) = \sum_{n = 1}^{\infty} \frac{z^n}{n^2} \end{equation}
	and can be extended to $\mathbb{C} \setminus (1; \infty)$ as
	\begin{equation}\label{li_2_integral} \Li(z) = -\int_{0}^{z} \frac{\ln(1-u)}{u} du, \end{equation}
	where the integral can be taken along the ray joining $0$ to $z$.
	
	\begin{figure}[h]
		\centering

		\begin{tikzpicture}
			\begin{axis}[
				axis lines = center,
				xlabel = $x$,
				ylabel = $y$,
				tick label style={/pgf/number format/fixed}
			]
				\addplot[
				    domain=0:1,
				    samples=100,
				    color=blue,
				    thick
				] {1/4 - (3 / (2 * pi^2)) * (x^2 + x^4 / 4 + x^6 / 9 + x^8 / 16 + x^10 / 25 + x^12 / 36)};
			\end{axis}
		\end{tikzpicture}

		\label{plot}
		\caption{Plot of $P(r)$}
	\end{figure}
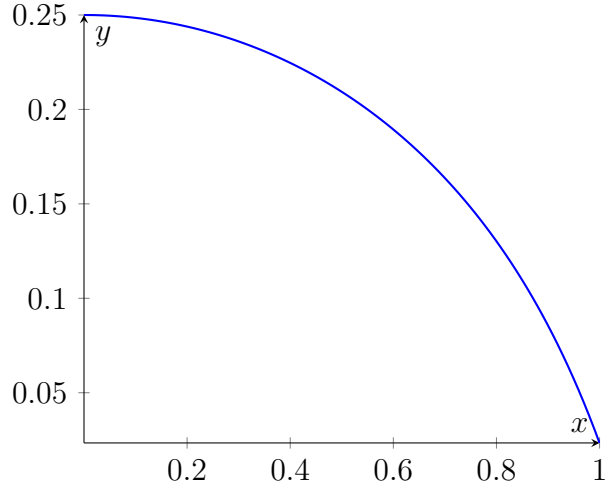

	Note that at $0$ and $1$ the dilogarithm function takes values
	\[ \Li(0) = 0, \quad \Li(1) = \sum_{n = 1}^{\infty} \frac{1}{n^2} = \frac{\pi^2}{6}, \]
	from which we see that (\ref{main_eq}) agrees with (\ref{end_values}).

	In this paper we present an alternative proof of Theorem \ref{main_theorem}. Through several manipulations we obtain formula (\ref{p_r_formula}) for the probability in question
	\[ P(r) = \frac{1}{4} - \frac{3}{4} \int_C \textrm{dif}(XC)^2 d\mu, \]
	where $\textrm{dif}(XC)$ is the difference of two arcs that the line $XC$ separates the circle into. The resulting integral expression we simply evaluate by means of contour integration.

	\textbf{Dilogarithm in geometric probability problems.} It is curious how the dilogarithm function arises as an answer to problems in geometric probability.

	For example, Karamata \cite{Karamata} has shown that the intersection point of two random chords of a unit circle lies within distance $r$ from the origin with probability
	\begin{equation}\label{chords} \frac{6}{\pi^2} \Li(r^2) \end{equation}
	for $r \in [0;1]$. This result was extended to $r > 1$ in \cite{Bortolotto}. Is there any connection between (\ref{main_eq}) and (\ref{chords})?

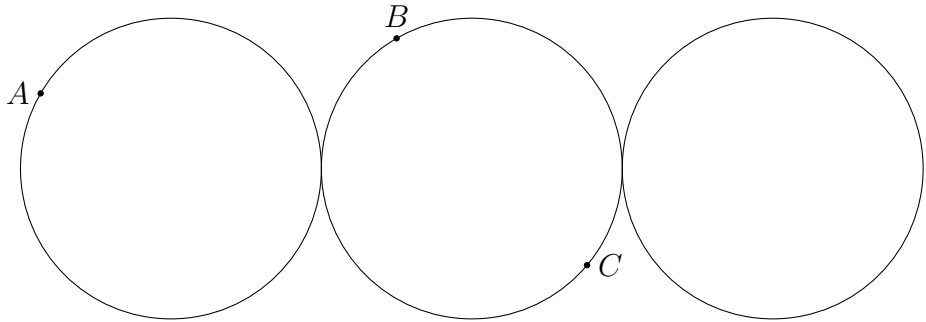
\begin{figure}[h]
		\centering
		\begin{tikzpicture}
			\draw (0, 0) circle (2);
			\draw (-4, 0) circle (2);
			\draw (4, 0) circle (2);

			\draw[fill] (+120: 2) circle (1pt) node[above] (B) {$B$};
			\draw[fill] (-40: 2) circle (1pt) node[right] (C) {$C$};

			\draw[fill, shift=({-4,0})] (150: 2) circle (1pt) node[left] (A) {$A$};
		\end{tikzpicture}
		\caption{Three tangent circles}
		\label{three_circles}
	\end{figure}

	A mysterious coincidence has been observed in a Mathematics Stack Exchange post.\footnote{\url{https://math.stackexchange.com/questions/5089940/why-do-these-two-lines-have-the-same-probability-of-intersecting-the-circle}} Take three circles of equal radii tangent to each other as in Figure \ref{three_circles}. Their centres are collinear and distinct.

	Random point $A$ is chosen on the leftmost circle. Random points $B$ and $C$ are chosen on the middle circle. Turns out that lines $AB$ and $BC$ have the same probability of intersecting the rightmost circle. The probability is
	\[ \frac{3}{4} - \frac{2}{\pi^2} \left( \frac{3}{8} \ln^2 2 + \Li(-\sqrt{2}) + 3 \Li\left( \frac{1}{\sqrt{2}} \right) \right) \]
	which is approximately $0.3871287106$.

	\subsection{Outline of the proof}

	Notation introduced in Section \ref{notation} would play a key role in the proof of Theorem \ref{main_theorem}. In Section \ref{condit} we use several non-trivial relations (\ref{relation_1}), (\ref{relation_2}), (\ref{relation_3}) between (\ref{not_1}), (\ref{not_2}), (\ref{not_3}) to derive a linear equation (\ref{lin_eq}) on $P(r)$. We then express $P(r)$ as a certain integral (\ref{p_r_exact}). The latter we evaluate in Section \ref{integ}.

	\subsection{Notation}
	\label{notation}
	Assume our circle has unit radius $R = 1$. Let $A, B$ and $C$ be points on the circle and let $X$ be a point inside the circle. By $\overrightarrow{AB}$ denote the line $AB$ oriented from $A$ to $B$.
	\begin{equation}\label{not_1} [X \in ABC] = \begin{cases} 1, & \textrm{if } X \textrm{ is inside triangle } ABC \\ 0, & \textrm{ otherwise}. \end{cases} \end{equation}
	\begin{equation}\label{not_2} \mathcal{R}(X, \overrightarrow{AB}) = \begin{cases} +1, & \textrm{if $X$ is to the right of $\overrightarrow{AB}$ or on $\overrightarrow{AB}$} \\ -1, & \textrm{if $X$ is to the left of $\overrightarrow{AB}$}. \end{cases} \end{equation}
	\begin{equation}\label{not_3} \widehat{ABC} = \begin{cases} +1, & \textrm{if triangle $ABC$ is oriented clockwise} \\ -1, & \textrm{if triangle $ABC$ is oriented counterclockwise}. \end{cases} \end{equation}

	An integral over the set of triangles we denote by
	\[ \int_{ABC} f(ABC) = \int_{A} \int_{B} \int_{C} f(ABC) (d\mu)^3,  \]
	where $\mu$ is the normalized Lebesgue measure on the unit circle. Also an integral over the set of triangles oriented clockwise we denote by
	\[ \int_{\widehat{ABC} = 1} f(ABC) = \iiint_{\widehat{ABC} = 1} f(ABC) (d\mu)^3. \]

	Let
	\[ \angle(\overrightarrow{MN}, \overrightarrow{KL}) \]
	be the lesser angle between a pair of vectors $\overrightarrow{MN}$ and $\overrightarrow{KL}$. The argument of a complex number $z \in \mathbb{C} \setminus (-\infty; 0]$ we denote by $\arg(z)$, on the positive real axis $\arg(z)$ is zero.

	\section{Proof of Theorem \ref{main_theorem}}
	\subsection{Calculus of conditions}
	\label{condit}

	Observe that $[X \in ABC]$ can be written in terms of $\mathcal{R}$. Indeed, let $ABC$ be a triangle oriented clockwise and consider sum $s = \mathcal{R}(X, \overrightarrow{AB}) + \mathcal{R}(X, \overrightarrow{BC}) + \mathcal{R}(X, \overrightarrow{CA})$. On the inside of $ABC$ it is equal $3$ and on the outside it is equal $1$. If we give the same triangle a counterclockwise orientation we would have $-3$ on the inside and $-1$ on the outside instead. We thus obtain a relation
	\begin{equation}\label{relation_1} (2 [X \in ABC] + 1) \widehat{ABC} = \mathcal{R}(X, \overrightarrow{AB}) + \mathcal{R}(X, \overrightarrow{BC}) + \mathcal{R}(X, \overrightarrow{CA}). \end{equation}
	\begin{figure}[h]
		\centering

		\begin{tikzpicture}
			\draw (0, 0) circle (60pt);
			\draw[fill] (160: 60pt) circle (1pt) node[left, scale=1] (A) {$A$};
			\draw[fill] (30: 60pt) circle (1pt) node[right, scale=1] (B) {$B$};
			\draw[fill] (-110: 60pt) circle (1pt) node[below, scale=1] (C) {$C$};

			\draw[-{Stealth[scale=1.5]}, thick] (160: 60pt) -- (30: 60pt);
			\draw[-{Stealth[scale=1.5]}, thick] (30: 60pt) -- (-110: 60pt);
			\draw[-{Stealth[scale=1.5]}, thick] (-110: 60pt) -- (160: 60pt);

			\node[scale=1.5] at (-0.3, -0.1) {$3$};
			\node[scale=1.5] at (0, 1.5) {$1$};
			\node[scale=1.5] at (1.1, -1) {$1$};
			\node[scale=1.5] at (-1.6, -0.8) {$1$};
		\end{tikzpicture}

		\label{sum_observation}
		\caption{}
	\end{figure}

	A map $ABC \to BAC$ that swaps a pair of vertices is a measure-preserving bijection between clockwise oriented triangles and counterclockwise oriented triangles. We can then restrict our attention to clockwise oriented triangles
	\begin{multline*}
		P(r) = \int_{ABC} [X \in ABC] = \int_{\widehat{ABC} = 1} [X \in ABC] + \int_{\widehat{ABC} = -1} [X \in ABC] \\ = \int_{\widehat{ABC} = 1} [X \in ABC] + \int_{\widehat{BAC} = 1} [X \in BAC] = 2\int_{\widehat{ABC} = 1} [X \in ABC]
	\end{multline*}
	and apply relation (\ref{relation_1})
	\begin{equation*}
		\int_{\widehat{ABC}=1} 2[X \in ABC] = \int_{\widehat{ABC}=1} [\mathcal{R}(X, \overrightarrow{AB}) + \mathcal{R}(X, \overrightarrow{BC}) + \mathcal{R}(X, \overrightarrow{CA}) - 1].
	\end{equation*}
	By symmetry
	\[ \int_{\widehat{ABC}=1} \mathcal{R}(X, \overrightarrow{AB}) = \int_{\widehat{ABC}=1} \mathcal{R}(X, \overrightarrow{BC}) = \int_{\widehat{ABC}=1} \mathcal{R}(X, \overrightarrow{CA}). \]
	In summary, we have
	\begin{equation}\label{R_integral} P(r) = 3 \int_{\widehat{ABC} = 1} [\mathcal{R}(X, \overrightarrow{AB})] - \frac{1}{2}. \end{equation}

	Note that a triangle $ABC$ is clockwise oriented if and only if $C$ is to the right of the chord $\overrightarrow{AB}$
	\begin{equation}\label{relation_2} \widehat{ABC} = \mathcal{R}(C, \overrightarrow{AB}). \end{equation}
	Using this relation we can extend the domain of integration in (\ref{R_integral}) to the set of all triangles
	\begin{equation}\label{domain_int} \int_{\widehat{ABC} = 1} \mathcal{R}(X, \overrightarrow{AB}) = \int_{ABC} \mathcal{R}(X, \overrightarrow{AB}) \frac{1 + \mathcal{R}(C, \overrightarrow{AB})}{2}. \end{equation}

	In an integral
	\begin{equation}\label{vanish_integral} \int_{ABC} \mathcal{R}(X, \overrightarrow{AB}) \end{equation}
	we can swap vertices $A$ and $B$. Since for a point $X$ not on $AB$
	\[ \mathcal{R}(X, \overrightarrow{BA}) = -\mathcal{R}(X, \overrightarrow{AB}), \]
	the integral (\ref{vanish_integral}) vanishes.

	Yet another relation would help us transform integral
	\[ \int_{ABC} \mathcal{R}(X, \overrightarrow{AB}) \mathcal{R}(C, \overrightarrow{AB}). \]
	Observe that products $\mathcal{R}(X, \overrightarrow{AB}) \mathcal{R}(C, \overrightarrow{AB})$ and $\mathcal{R}(A, \overrightarrow{XC}) \mathcal{R}(B, \overrightarrow{XC})$ yield an expression for $[X \in ABC]$ (assuming $X$ does not lie on triangle $ABC$)
	\begin{equation}\label{relation_3} 2[X \in ABC] = \mathcal{R}(X, \overrightarrow{AB}) \mathcal{R}(C, \overrightarrow{AB}) - \mathcal{R}(A, \overrightarrow{XC}) \mathcal{R}(B, \overrightarrow{XC}), \end{equation}
	which we verify by considering two cases: segments $AB$ and $CX$ intersect or segments $AB$ and $CX$ do not intersect. Assume we are in the first case. Points $C$ and $X$ are on different sides of $AB$, so $X$ cannot be inside $ABC$. Also points $A$ and $B$ are on different sides of $CX$. In other words, both products $\mathcal{R}(X, \overrightarrow{AB}) \mathcal{R}(C, \overrightarrow{AB})$ and $\mathcal{R}(A, \overrightarrow{XC}) \mathcal{R}(B, \overrightarrow{XC})$ are equal $-1$, and equality in (\ref{relation_3}) follows. Now assume we are in the second case. Points $C$ and $X$ are on the same side of $AB$, meaning that $\mathcal{R}(X, \overrightarrow{AB}) \mathcal{R}(C, \overrightarrow{AB})$ is equal to $1$. Point $X$ is inside triangle $ABC$ if and only if $X$ is inside angle $\angle BCA$, i.e., line $CX$ separates $A$ and $B$. The latter is equivalent to $\mathcal{R}(A, \overrightarrow{XC}) \mathcal{R}(B, \overrightarrow{XC})$ being equal $-1$, and equality in (\ref{relation_3}) again follows.
	\begin{figure}[h]
		\centering		
		\begin{subfigure}{0.49\textwidth}
			\centering

			\begin{tikzpicture}
				\draw (0, 0) circle (60pt);
				\draw[fill] (-60: 60pt) circle (1pt) node[right] (A) {$A$};
				\draw[fill] (80: 60pt) circle (1pt) node[above] (B) {$B$};
				\draw[fill] (-170: 60pt) circle (1pt) node[left] (C) {$C$};
				\draw[fill] (30: 40pt) circle (1pt) node[right] (X) {$X$};

				\draw (-60: 60pt) -- (80: 60pt) -- (-170: 60pt) -- (-60: 60pt);
				\draw (-170: 60pt) -- (30: 40pt);
			\end{tikzpicture}
			\caption{Segments $AB$ and $CX$ intersect}
			\label{ab_cx_cases:first}
		\end{subfigure}
		\hfill
		\begin{subfigure}{0.49\textwidth}
			\centering

			\begin{tikzpicture}
				\draw (0, 0) circle (60pt);
				\draw[fill] (40: 60pt) circle (1pt) node[right] (A) {$A$};
				\draw[fill] (140: 60pt) circle (1pt) node[above] (B) {$B$};
				\draw[fill] (-170: 60pt) circle (1pt) node[left] (C) {$C$};
				\draw[fill] (-80: 40pt) circle (1pt) node[right] (X) {$X$};

				\draw (40: 60pt) -- (140: 60pt) -- (-170: 60pt) -- (40: 60pt);
				\draw (-170: 60pt) -- (-80: 40pt);
			\end{tikzpicture}
			\caption{Segments $AB$ and $CX$ do not intersect}
			\label{ab_cx_cases:second}
		\end{subfigure}
		\caption{}
	\end{figure}

 	Applying (\ref{relation_3}) we get
 	\begin{multline}\label{integ_sq} \int_{ABC} \mathcal{R}(X, \overrightarrow{AB}) \mathcal{R}(C, \overrightarrow{AB}) = \int_{ABC} [2[X \in ABC] + \mathcal{R}(A, \overrightarrow{XC}) \mathcal{R}(B, \overrightarrow{XC})] \\ = 2 P(r) + \int_{ABC} \mathcal{R}(A, \overrightarrow{XC}) \mathcal{R}(B, \overrightarrow{XC}) \\ = 2P(r) + \int_{C} \left( \int_{A} \mathcal{R}(A, \overrightarrow{XC}) d\mu \right)^2 d\mu, \end{multline}
 	where the last equality holds, since $\mathcal{R}(A, \overrightarrow{XC})$ and $\mathcal{R}(B, \overrightarrow{XC})$ are two independent functions in variables $A$ and $B$.

 	Putting (\ref{R_integral}), (\ref{domain_int}), (\ref{integ_sq}) together we obtain
 	\begin{equation}
 		\label{lin_eq}
 		P(r) = -\frac{1}{2} + \frac{3}{2} \int_{C} \left( \int_{A} \mathcal{R}(A, \overrightarrow{XC}) d\mu \right)^2 d\mu + 3P(r).
 	\end{equation}
 	Solving this linear equation we get a formula for $P(r)$
 	\begin{equation}
 		\label{p_r_formula}
 		P(r) = \frac{1}{4} - \frac{3}{4} \int_{C} \left( \int_{A} \mathcal{R}(A, \overrightarrow{XC}) d\mu \right)^2 d\mu.
 	\end{equation}

	\subsection{The integral}
	\label{integ}

	Line $XC$ cuts our circle into two arcs. On one arc $\mathcal{R}(A, \overrightarrow{XC}) = 1$ and on another $\mathcal{R}(A, \overrightarrow{XC}) = -1$. The difference of two arcs is equal $2\angle MON$ (see Figure \ref{chase}). The central angle $\angle MON$ is double the measure of the inscribed angle $\angle OCX$. Thus
 	\[ \int_{A} \mathcal{R}(A, \overrightarrow{XC}) d\mu = \pm \frac{1}{2\pi} (4 \angle OCX) = \pm \frac{2}{\pi} \angle (\overrightarrow{OC}, \overrightarrow{XC}). \]

 	\begin{figure}[h]
 		\centering	
 		\begin{tikzpicture}
 			\draw (0, 0) circle (60pt);
			\draw[fill] (0, 0) circle (1pt) node[right] (O) {$O$};
			\draw[fill] (0 : 30pt) circle (1pt) node[right] (X) {$X$};
			\draw[fill] (50 : 60pt) circle (1pt) node[right] (A) {$C$};
			\draw[fill] (-71 : 60pt) circle (1pt) node[below] (N) {$N$};
			\draw[fill] (-130 : 60pt) circle (1pt) node[below left] (M) {$M$};

			\draw (50 : 60pt) -- (-71 : 60pt);
			\draw (0, 0) -- (50 : 60pt);
			\draw (50 : 60pt) -- (-130 : 60pt);
			\draw (0, 0) -- (-71 : 60pt);
 		\end{tikzpicture}

 		\caption{}
 		\label{chase}
 	\end{figure}

 	Now we switch to complex numbers. Let $\overrightarrow{OC}$ correspond to a number $e^{i \theta}$ on the unit circle in the complex plane and let $\overrightarrow{OX}$ correspond to $r$. In this new notation,
 	\begin{multline*} \angle (\overrightarrow{OC}, \overrightarrow{XC}) = \angle(c, c - r) = \angle(1, 1 - r/c) \\ = \arg(1 - r e^{-i \theta}) = -\arg(1 - r e^{i \theta}). \end{multline*}
 	In the second equality we use that a rotation by $-\theta$ leaves the angle invariant, and in the last equality we use conjugation
 	\[ \arg(\overline{z}) = -\arg(z). \]
 	Formula (\ref{p_r_formula}) becomes
 	\begin{equation}\label{p_r_exact} P(r) = \frac{1}{4} - \frac{3}{2\pi^3} \int_{0}^{2\pi} \arg(1 - r e^{i\theta})^2 d\theta. \end{equation}

 	The integral in question
 	\[ I(r) = \int_{0}^{2 \pi} \arg(1 - r e^{i\theta})^2 d \theta \]
 	we differentiate at $r \in (0;1)$ to get rid of the square
 	\begin{equation}\label{d_dr} \frac{d}{d r} I(r) = -2 \int_{0}^{2 \pi} \frac{\sin \theta}{1 + r^2 - 2r \cos \theta} \arg(1 - r e^{i \theta}) d\theta. \end{equation}
 	To elaborate, note that $1 - r e^{i \theta}$ lies in the open right half-plane $\operatorname{Re} z > 0$, so $\arg(1 - r e^{i \theta})$ is the imaginary part of $\ln(1 - re^{i \theta})$, where the latter is the complex logarithm function with domain $\mathbb{C} \setminus (-\infty; 0]$ and $\ln(1) = 0$. This observation helps us evaluate the derivative of $\arg(1 - r e^{i\theta})$
 	\begin{multline*} \frac{\partial}{\partial r} \arg(1 - re^{i\theta}) = \operatorname{Im}\left( \frac{\partial}{\partial r} \ln(1 - re^{i\theta}) \right) = \operatorname{Im}\left( \frac{-e^{i\theta}}{1 - r e^{i \theta}} \right) \\ = \operatorname{Im}\left( \frac{r - e^{i\theta}}{(1 - r e^{i \theta})(1 - r e^{-i\theta})} \right) = \frac{-\sin \theta}{1 + r^2 - 2r \cos \theta}. \end{multline*}

 	Next, we turn (\ref{d_dr}) into a contour integral
 	\begin{multline*} \frac{d}{dr} I(r) = -2 \int_{0}^{2 \pi} \frac{\frac{1}{2 i}(e^{i\theta} - e^{-i\theta})}{(1 - r e^{i\theta})(1 - r e^{-i\theta}) (i e^{i \theta})} \arg(1 - r e^{i \theta}) (i e^{i \theta} d\theta) \\ = \int_{0}^{2\pi} \frac{e^{i\theta} - e^{-i \theta}}{(1 - re^{i\theta})(e^{i\theta} - r)} \arg(1 - r e^{i\theta}) (i e^{i\theta} d\theta) \\ = \int_{0}^{2\pi} \frac{e^{2 i \theta} - 1}{(1 - re^{i\theta})(e^{i\theta} - r) e^{i\theta}} \arg(1 - r e^{i\theta}) (i e^{i\theta} d\theta) \\ = \operatorname{Im} \oint_{|z| = 1} \frac{z^2 - 1}{(1-rz)(z-r)z} \ln(1 - rz) dz. \end{multline*}
 	The complex logarithm function is analytic in the open right half-plane, so we can use the residue theorem. Function
 	\[ \frac{z^2 - 1}{(1-rz)(z-r)z} \ln(1 - rz) \]
 	has poles $z = 0$ and $z = r$ within the contour of integration. At $z = 0$ the residue vanishes, since $\ln(1) = 0$, and at $z = r$ the residue is equal to
 	\[ \frac{r^2 - 1}{(1 - r^2) r} \ln(1 - r^2) = -\frac{\ln(1-r^2)}{r}. \]
 	Thus our contour integral evaluates to
 	\[ \frac{d}{dr} I(r) = - 2 \pi \frac{\ln(1-r^2)}{r}. \]

 	To restore $I(r)$ we integrate the above equality
 	\begin{multline*} I(r) = I(0) - 2 \pi \int_{0}^{r} \frac{\ln(1 - r^2)}{r} dr = -\pi \int_{0}^{r} \frac{\ln(1-r^2)}{r^2} (2rdr) \\ = -\pi \int_{0}^{r^2} \frac{\ln(1 - u)}{u} du = \pi \Li(r^2). \end{multline*}
 	Substituting this into (\ref{p_r_exact}) we arrive at
 	\[ P(r) = \frac{1}{4} - \frac{3}{2 \pi^2} \Li(r^2). \]

\end{document}